\documentclass[a4paper,12pt]{article}
\usepackage{amsmath}
\usepackage{amssymb}
\usepackage{amsfonts}
\begin{document}
\newtheorem{theoremenv}{Theorem}
\newtheorem{propositionenv}[theoremenv]{Proposition}
\newtheorem{lemmaenv}[theoremenv]{Lemma}
\newtheorem{corollaryenv}[theoremenv]{Corollary}
\newtheorem{claimenv}[theoremenv]{Claim}
\newtheorem{constructionenv}[theoremenv]{Construction}
\newtheorem{conjectureenv}[theoremenv]{Conjecture}
\newtheorem{questionenv}[theoremenv]{Question}
\newtheorem{definitionenv}[theoremenv]{Definition}
\newtheorem{problemenv}[theoremenv]{Problem}
\newsavebox{\proofbox}
\savebox{\proofbox}{\begin{picture}(7,7)%
  \put(0,0){\framebox(7,7){}}\end{picture}}

\def\proof{\noindent\textbf{Proof. }}
\def\endproof{\hfill{\usebox{\proofbox}}}
\newcommand{\la}{\left\langle}
\newcommand{\lv}{\left\Vert}
\newcommand{\lmod}{\left|}
\newcommand{\lb}{\left(}
\newcommand{\lc}{\left\{}
\newcommand{\ra}{\right\rangle}
\newcommand{\rv}{\right\Vert}
\newcommand{\rmod}{\right|}
\newcommand{\rb}{\right)}
\newcommand{\rc}{\right\}}
\newcommand{\Supp}{\mbox{Supp}}
\begin{center}
\Large On the Hardy-Littlewood majorant problem\\[15pt]
\large Ben Green\footnote{While this work was carried out the first author was resident in Budapest, and was supported by the \textit{Mathematics in Information Society} project carried out by R\'enyi Institute, in the framework of the European Community's \textit{Confirming the International R\^ole of Community Research} programme.} and Imre Z. Ruzsa\\[20pt]
\end{center}
\begin{abstract}\noindent Let $\Lambda \subseteq \{1,\dots,N\}$, and let $\{a_n\}_{n \in \Lambda}$ be a sequence with $|a_n| \leq 1$ for all $n$. It is easy to see that 
\[ \left\Vert \sum_{n \in \Lambda} a_n e(n\theta) \right\Vert_{p} \; \leq \; \left\Vert \sum_{n \in \Lambda} e(n\theta) \right\Vert_{p}\] for every even integer $p$. We give an example which shows that this statement can fail rather dramatically when $p$ is not an even integer. This answers in the negative a question known as the Hardy-Littlewood majorant conjecture, thereby ruling out a certain approach to the restriction and Kakeya families of conjectures.\end{abstract}
\normalsize \noindent\textbf{1. Introduction} Let $\Lambda \subseteq \{1,\dots,N\}$, and let $\{a_n\}_{n \in \Lambda}$ be a sequence with $|a_n| \leq 1$ for all $n$. Hardy and Littlewood observed (as a simple consequence of Parseval's identity) that we have \begin{equation}\label{star1} \left\Vert \sum_{n \in \Lambda} a_n e(n\theta) \right\Vert_{p} \; \leq \; \left\Vert \sum_{n \in \Lambda} e(n\theta) \right\Vert_{p}\end{equation} for every even integer $p$. They asked what, if anything, can be said about other values of $p > 2$. Defining a constant $B_p(\Lambda)$ by 
\begin{equation}\label{star2} \sup_{\{a_n\} \, : \, |a_n| \leq 1}\left\Vert \sum_{n \in \Lambda} a_n e(n\theta) \right\Vert_{p} \; = \; B_p(\Lambda)\left\Vert \sum_{n \in \Lambda} e(n\theta) \right\Vert_{p},\end{equation} their question may be interpreted as asking for the behaviour of $B_p(\Lambda)$. Thus for any $\Lambda$ one has $B_p(\Lambda) = 1$ for all even integers $p$. As Hardy and Littlewood \cite{HL1,HL2} knew, it is possible for $B_3(\Lambda)$ to be larger than $1$, so that a perfect analogue of \eqref{star1} cannot hold. We will discuss an example later on.\\[11pt]
Let us give a brief history of the problem since Hardy and Littlewood. It is natural to write
\[ B_p(N) \; = \; \sup_{\Lambda \subseteq \{1,\dots,N\}} B_p(\Lambda)\] and to ask for the behaviour of $B_p(N)$. Thus Hardy and Littlewood knew that $B_3(N) > 1$, and
Boas \cite{Boas} later showed that $B_p(N) > 1$ for \textit{any} $p \notin \{2,4,6,\dots\}$. Disproving a conjecture of Hardy and Littlewood, Bachelis \cite{Bachelis} (see also \cite{Mont} p. 138) showed that in fact $B_p(N) \rightarrow \infty$ for any $p \notin \{2,4,6,\dots\}$. The idea behind this construction is to take a kind of product of Boas-type examples, an idea which Bachelis states was communicated to him by Katzelson. A more precise version of the same idea forms the heart of the present paper.\\[11pt]
Recently, there has been renewed interest in the quantitative behaviour of $B_p(N)$ due to connections with two famous open questions in harmonic analysis, the restriction conjecture and the Kakeya conjecture. The Kakeya conjecture is a problem in geometric measure theory of much greater importance than one would at first sight think. In its weakest form, it asserts that a Besicovitch set, that is any compact subset of $\mathbb{R}^d$ containing a unit line segment in each direction, must have Minkowski dimension $d$. The restriction conjecture is, put simply, an assertion about the Fourier transforms of functions supported on the $(d-1)$-dimensional Euclidean sphere $S^{d-1} \subseteq \mathbb{R}^d$. One form of it is the following.
\begin{conjectureenv}[Local restriction conjecture]\label{conj1} There exists, for each $\epsilon > 0$,  an absolute constant $\gamma = \gamma_{d,\epsilon}$ with the following property. Let $\sigma$ be the surface measure on $S^{d-1}$, and suppose that $f : S^{d-1} \rightarrow \mathbb{C}$ is a continuous function with $\Vert f \Vert_{\infty} \leq 1$. Then, for all $R$ sufficiently large,
\[ \int_{|\xi| \leq R} |\widehat{f d\sigma}(\xi)|^{2d/(d-1)}\, d\xi \; \leq \; \gamma R^{\epsilon}.\] 
\end{conjectureenv}
This implies a rather more classical global restriction conjecture of the form 
\begin{equation}\label{zz} \int |\widehat{f d\sigma}(\xi)|^{2d/(d-1) + \delta} \, d\xi \; = \; O_{\delta,d}(1)\end{equation} by the so-called $\epsilon$-removal technique of Tao \cite{ResBoc}.
Conjecture \ref{conj1} is known to be true for $d = 2$ (see \cite{Zyg}) but not for any higher value of $d$. It would, by a classical method (essentially due to Fefferman -- see \cite{GreenRKP}, Chapter 5) imply the Kakeya conjecture.\\[11pt]
The novice reader may be rather mystified by all of this, but by now there are several good introductions to this area of research, such as \cite{Tao1,Tao2}. See also \cite{GreenRKP}.\\[11pt]
By an argument which is sketched in \cite{Mock}, Conjecture \ref{conj1} would follow if we could show that \begin{equation}\label{hopeful}B_p(\Gamma) \; \ll_{p,\epsilon} \; N^{\epsilon},\end{equation} where $\Gamma \subseteq \{1,\dots,N\}$ is a set constructed by discretizing and then projecting the sphere $S^{d-1}$ in a suitable way. The resulting set $\Gamma$ is not particularly natural, but this would not concern us if we could prove that in fact \eqref{hopeful} holds for \textit{any} $\Lambda$ in place of $\Gamma$, that is that $B_p(N) \ll_{p,\epsilon} N^{\epsilon}$. The question of whether such a statement holds has become known as the Hardy-Littlewood majorant problem and, as we have stated, it would imply the restriction conjecture, and hence the Kakeya conjecture.
\begin{problemenv}[Hardy-Littlewood majorant problem]\label{mainprob} Is it true that we have an estimate
\[ B_p(N) \; \ll_{p,\epsilon} \; N^{\epsilon}\] for all $\epsilon > 0$ and $p \geq 2$?
\end{problemenv}
In \cite{Mock} it is shown that $B_p(N) \gg \exp(c\log N/\log\log N)$. In \cite{MockSchlag} a detailed investigation of the behaviour of $B_p(\Lambda)$ for \textit{random} subsets $\Lambda \subseteq \{1,\dots,N\}$ is conducted. In particular it is shown that for any fixed $\lambda \in (0,1)$ and for any $p \in (2,4)$ (the range of interest for applications) a random subset $\Lambda \subseteq \{1,\dots,N\}$ with cardinality $N^{\lambda}$ will satisfy $B_p(\Lambda) \ll_{p,\epsilon,\lambda} N^{\epsilon}$ almost surely.\\[11pt]
The main theorem of this paper is that the answer to Problem \ref{mainprob} is \textit{no}. 
\begin{theoremenv}\label{mainthm} $B_3(N) \gg N^{\eta}$ for some explicitly computable positive constant $\eta$.\end{theoremenv}
The proof is by construction of a set $\Lambda$ with $B_3(\Lambda)$ large, and is completely explicit. In the next section we describe the construction, and then we offer two alternative proofs of Theorem \ref{mainthm}. \\[11pt]
We have recently heard that Mockenhaupt and Schlag have independently obtained a proof of Theorem \ref{mainthm}. In fact, they obtain the same result with $3$ replaced by any $p > 2$ which is not an even integer. Our construction could be modified to obtain such a result, but we have not done so here. The two approaches are quite similar, which is perhaps rather unsurprising. Mockenhaupt and Schlag's argument will appear in \cite{MockSchlag}.\\[11pt]
\noindent\textbf{2. The construction.} It seems to be quite well-known (see \cite{Mont}, p144) that if $Q(\theta) = 1 + e(\theta) + e(3\theta)$ and $q(\theta) = 1 + e(\theta) - e(3\theta)$ then
\[ \int^1_0 |q(\theta)|^3\,d\theta \; > \; \int^1_0  |Q(\theta)|^3\,d\theta.\] In any case this statement is rather easy to check on a computer, and it confirms that the set $\Lambda = \{0,1,3\}$ has $B_3(\Lambda) > 1$. It is rather natural to try taking a product of several copies of this set. Thus if $D$ is a large positive integer, define the set
\[ \Lambda \; = \; \left\{ \sum_{i = 0}^{k-1} \varepsilon_i D^i \; | \; \mbox{$\varepsilon_i \in \{0,1,3\}$ for all $i$}\right\},\] where $k = \lfloor \log N/\log D\rfloor$. We will show that if $D$ is large enough then these sets satisfy Theorem \ref{mainthm}. In fact if $a_n$ is defined to be $(-1)^{W(n)}$, where $W(n)$ is the number of 3s in the base $D$ expansion of $D$, then we will be able to prove that $\Vert \sum a_n e(n\theta) \Vert_3$ is much larger than $\Vert \sum e(n\theta) \Vert_3$. It is easy to see that
\begin{equation}\label{eq209} \sum_{n \in \Lambda} e(n\theta) \; = \; \prod_{i=0}^{k-1} Q(D^i \theta)\end{equation} and
\begin{equation}\label{eq210} \sum_{n \in \Lambda} a_n e(n\theta) \; = \; \prod_{i=0}^{k-1} q(D^i \theta).\end{equation}
We we will write these two expressions as $F(\theta)$ and $f(\theta)$ respectively.\\[11pt]
\noindent\textbf{3. Proof that the construction works.}
To estimate the norms of $F$ and $f$ we prove the following.
\begin{lemmaenv}
Let $p$ be \emph{any} trigonometric polynomial and, for let $D$ be a large positive integer. Define a trigonometric polynomial $g$ by
\[ g(\theta) \; = \; \prod_{i=0}^{k-1} p(D^i \theta)\]
\emph{(}viz. equation \eqref{eq209}\emph{)}. Let $\alpha$ and $\varepsilon$ be positive numbers. Then there
is a $D_0 = D_0(p,\alpha,\varepsilon)$ such that for all $D>D_0$ we have
\begin{equation} \label{a}
\Vert p\Vert _\alpha  - \varepsilon  \; \leq \; \Vert g\Vert _\alpha ^{1/k} \; \leq \; \Vert p\Vert _\alpha  + \varepsilon   
\end{equation}
for every $k$.
\end{lemmaenv}
\noindent\textbf{Remark.} Of course, we are interested in the case $p = 1 + e(\vartheta) + e(3\vartheta)$ or $1 + e(\theta) - e(3\vartheta)$.\\[11pt]
\proof 
We proceed by finding functions $r^+$ and $r^-$ satisfying
\begin{equation} \label{b}
0  \; \leq  \; r^-(\phi) \; \leq \;  |p(\phi)|^\alpha \;  \leq \; r^+(\phi) 
\end{equation}
and
\begin{equation} \label{c}
 (\Vert p\Vert _\alpha  - \varepsilon )^\alpha \; \leq \; \Vert r^-\Vert _1 \; \leq \; \Vert r^+\Vert _1 \; \leq \; (\Vert p\Vert _\alpha  + \varepsilon )^\alpha 
\end{equation}
with the property that both $r=r^+$ and $r=r^-$ satisfy
\begin{equation} \label{d}
     \int _0^1 \prod _{i=0}^{k-1} r(D^i\vartheta ) d\vartheta  = \left ( \int _0^1 r(\vartheta ) d\vartheta  \right )^k 
     \end{equation} for any positive integer $k$.
Once such functions have been given we can estimate $\Vert g\Vert_{\alpha}$ by
     \begin{eqnarray*}  \Vert g\Vert_\alpha^{\alpha} & = & 
     \int _0^1 \prod _{i=0}^{k-1} |p(D^i\vartheta)|^{\alpha} d\vartheta \\ & \leq &
     \int _0^1 \prod _{i=0}^{k-1} r^+(D^i\vartheta ) d\vartheta \\ & =  & \left ( \int _0^1 r^+(\vartheta ) d\vartheta  \right )^k \\ &
     \leq & (\Vert p\Vert _\alpha  + \varepsilon )^{\alpha k},   \end{eqnarray*}
     which is the upper bound of \eqref{a}. The lower bound follows similarly by
using $r^-$.\\[11pt]
Properties \eqref{b} and \eqref{c} will follow from an assumption
     \begin{equation} \label{e}
     \max (0, |p(\phi)|^\alpha -\delta ) \; \leq \; r^-(\phi) \; \leq \; |p(\phi)|^\alpha  \; \leq \; r^+(\phi) \; \leq \; |p(\phi)|^\alpha  +\delta  
     \end{equation}
     with a suitably chosen positive $\delta $.\\[11pt]
One class of functions $r$ that have property \eqref{d} is $\mathcal{A}$, the functions that are
constant on each interval $R_j = [j/D, (j+1)/D)$. To see this, observe that \begin{equation}\label{gg} (\vartheta,D\vartheta,\dots,D^{k-1}\vartheta) \; \in \; R_{j_0} \times \dots \times R_{j_{k-1}}\end{equation} precisely for those $\vartheta$ whose base $D$ expansion starts
\[ \vartheta \; = \; \frac{j_0}{D} + \frac{j_1}{D^2} + \dots + \frac{j_{k-1}}{D^k} + \dots\] Thus for fixed $j_0,\dots,j_{k-1}$ the inclusion \eqref{gg} holds for $\vartheta$ lying in some interval of length exactly $D^{-k}$, and this quickly implies \eqref{d}. The continuity of $p$ ensures that we may find $r \in \mathcal{A}$ which satisfy \eqref{e} as well, provided of course that $D$ is sufficiently large.\\[11pt]
Another class of functions $r$ that can be used is $\mathcal{B}$, the trigonometrical polynomials of degree less than $D/2$. It is easy to see that any $r \in \mathcal{B}$ satisfies \eqref{d}. If $p$ has no zero on the unit circle (which is the case for the polynomials which actually interest us) then \eqref{e} is an immediate consequence of Weierstrass' approximation theorem. 
If $p$ does have a zero, some care is needed to find the lower function $r^-$. One can circumvent the obstacle by 
seeking $r^-$ in the form
     \[   r^-(\phi) = |p(\phi)|^m h(\phi) ,  \]
     where $m > \alpha $ is an even integer and $h$ is a trigonometric polynomial. If $h$ satisfies
          \[  \max\left(0,\frac{|p(\phi)|^{\alpha}  - \delta}{|p(\phi)|^ m}\right) \; \leq \; h(\phi)\; \leq  \; \frac{ |p(\phi)|^\alpha  + \delta}{|p(\phi)|^ m} \] then $r^{-}$ will satisfy \eqref{e}.
     The existence of such an $h$ does now follow from Weierstrass'
theorem. \\[11pt]
Whether we use $\mathcal{A}$ or $\mathcal{B}$, the lemma has certainly been verified.\endproof\\[11pt]
We remark that the use of $\mathcal{A}$ is both
simpler and works for every continuous function $p$, while the use of $\mathcal{B}$ only works for functions that have only a finite number of roots, and if those roots have
finite order. It has, however, some hidden advantages. If $p$ is a nice
function (say, analytic), then approximation by polynomials of degree $D$ is
typically much better than by functions of step $1/D$. Furthermore, if we use
polynomials of higher degree, then, though the simple equality \eqref{d} will no
longer hold, it is still possible to calculate the integral on the left side by
a linear recurrence. This was done in a different context by Pintz and the second author
\cite{r03e}.\\[11pt]
By the remarks in \S 2, the deduction of Theorem \ref{mainthm} is simply a matter of proving a bound of the form $\Vert f \Vert_3 \gg N^{\eta}\Vert F \Vert_3$.
We apply the Lemma for the functions $p$ and $q$
defined in \S 2, $\alpha =3$ and 
     \[   \varepsilon  = ( \Vert q\Vert _3 - \Vert Q\Vert _3)/3 ,  \]
     say. Then we have $ \Vert f\Vert _3/\Vert F\Vert _3 \geq  c^k$ with
     \[   c = \frac{\Vert q\Vert _3 - \varepsilon}{ \Vert Q\Vert _3 + \varepsilon  } >1 .  \]
     Observing that $N<D^k$ we see that Theorem \ref{mainthm} holds with $\eta  = \log c /\log D$.\endproof\\[11pt]
     An explicit constant $\eta$ could be calculated, but it would be very small. To finish this section we show that $\eta \leq 1/18$ by modifying an argument sketched in \cite{Mock}.
\begin{propositionenv} Let $\Lambda \subseteq \{1,\dots,N\}$, and let $(a_n)_{n \in \Lambda}$ be a sequence with $|a_n| \leq 1$ for all $n$. Then
\[ \left\Vert \sum_{n \in \Lambda} a_n e(n\theta) \right\Vert_3 \; \ll \; N^{1/18}\left\Vert \sum_{n \in \Lambda} e(n\theta) \right\Vert_3.\]
\end{propositionenv}
\proof Consider the map $T : l^{\infty}(\Lambda) \rightarrow C(\mathbb{T})$ defined by
\[ (Ta)(\theta) \; = \; \sum a_n e(n\theta)\] for any $a = (a_n)_{n \in \Lambda}$. Write $F(\theta) = \sum_{n \in \Lambda} e(n\theta)$. Then it is easy to check, using Parseval's identity, the bounds
\[ \Vert Ta \Vert_2 \; \leq \; \Vert a \Vert_{\infty} \Vert F \Vert_2\]
and 
\[ \Vert Ta \Vert_4 \; \leq \; \Vert a \Vert_{\infty} \Vert F \Vert_4\] for any sequence $a$.
By the Riesz--Thorin interpolation theorem (one reference is \cite{GreenRKP}, Ch. 7.) these imply a bound
\begin{equation} \label{l3} \Vert Ta \Vert_3 \; \leq \; \Vert a \Vert_{\infty} \Vert F \Vert_2^{1/3} \Vert F \Vert_4^{2/3}.\end{equation}
Now observe that 
\begin{equation}\label{eq11} \Vert F \Vert_2 \; \leq \; \Vert F \Vert_3\end{equation}
and that 
\begin{eqnarray}\nonumber \Vert F \Vert_4 & \leq & \Vert F \Vert_{\infty}^{1/4} \Vert F \Vert_{3}^{3/4} \\ & \leq & |\Lambda|^{1/4} \Vert F \Vert_3^{3/4} \label{eq12}.\end{eqnarray}
Equations \eqref{l3},\eqref{eq11} and \eqref{eq12} immediately lead to
\begin{equation}\label{l4} \Vert Ta \Vert_3 \; \leq \; \Vert a \Vert_{\infty} |\Lambda|^{1/6}\Vert F \Vert_3^{5/6}.\end{equation}
However, observing that $|F(\theta)| \geq |\Lambda|/2$ for $|\theta| \leq 1/10N$, we see that $\Vert F \Vert_3 \gg |\Lambda|N^{-1/3}$. This provides an upper bound for $|\Lambda|$ which, when substituted into \eqref{l4}, gives the bound
\[ \Vert Ta \Vert_3 \; \ll \; N^{1/18}\Vert a \Vert_{\infty} \Vert F \Vert_3.\] This is what we wanted to prove.\endproof\\[11pt]
We remark if one wished to replace the exponent $3$ by some $p \in (2,4)$ then an identical argument would allow one to replace $1/18$ by $2(1/p - 1/4)(1 - 2/p)$.\\[11pt] 
\noindent\textbf{4. Further remarks on the majorant problem.} Our result, Theorem \ref{mainthm}, is rather negative since it does nothing more than rule out an approach to some well-known problems which was rather optimistic. The constants $B_p(\Lambda)$ are really of interest for specific sets which arise ``naturally'', for example from surfaces in Euclidean space.\\[11pt] Mockenhaupt \cite{Mock} proves that if $f : S^1 \rightarrow \mathbb{C}$ satisfies $\Vert f \Vert_{\infty} \leq 1$ then one has the inequality
\begin{equation}\label{eq200} \int_{|\xi| \leq R} |\widehat{fd\sigma}(\xi)|^4 \, d\xi \; \leq \; \pi \int_{|\xi| \leq R} |\widehat{d\sigma}(\xi)|^4 \, d\xi\end{equation} for $R$ sufficiently large (recall that $\sigma$ is the induced measure on $S^1$, considered as a subset of $\mathbb{R}^2$). This may be regarded as an estimate of Hardy-Littlewood majorant type, with the r\^ole of the set $\Lambda$ being played by $S^1$. We are not aware of any counterexamples to the corresponding statement in $d$ dimensions. That is, we cannot be certain that there is no absolute constant $C_d$ such that
\begin{equation}\label{eq20} \int_{|\xi| \leq R} |\widehat{fd\sigma}(\xi)|^{2d/(d-1)} \, d\xi \; \leq \; C_d \int_{|\xi| \leq R} |\widehat{d\sigma}(\xi)|^{2d/(d-1)} \, d\xi\end{equation} for all functions $f : S^{d-1} \rightarrow \mathbb{C}$ with $\Vert f \Vert_{\infty} \leq 1$ and all $R \geq R_0(d)$.\\[11pt]
The relevance of the exponent $2d/(d-1)$ is that $\widehat{d\sigma}$ lies in $L^q(\mathbb{R}^d)$ if and only if $q > 2d/(d-1)$.\\[11pt]
The estimate \eqref{eq20}, if true, would be quite remarkable. It would immediately imply a strong quantitative form of the local restriction conjecture (Conjecture \ref{conj1}). It would also imply, by the classical method that if $E \subseteq \mathbb{R}^d$ is a Besicovitch set then the $\delta$-neighbourhood $N_{\delta}(E)$ has measure at least $c_d\left(\log (1/\delta)\right)^{1-d}$. This would be an extremely strong quantitative form of the Kakeya conjecture. Trying to disprove this might not be a bad way of showing that \eqref{eq20} is too optimistic, if indeed this is the case.

\end{document}